\documentclass[12pt,reqno]{article}
\usepackage{amssymb}
\usepackage{graphicx}
\usepackage{amsmath}
\usepackage[latin1]{inputenc}
\usepackage{amssymb,indentfirst}
\usepackage{float}
\usepackage{times}
\usepackage[T1]{fontenc}
\usepackage{tikz}
\newtheorem{thm}{Theorem}[section]

\newtheorem{rem}[thm]{Remark}

\renewcommand{\raggedright}{\leftskip=0pt \rightskip=0pt plus 0cm}
\raggedright

\def\N{{\mathbb N}}

\def\R{{\mathbb R}}

\def\bb{\begin}
\def\bc{\begin{center}}
\def\ec{\end{center}}
\def\be{\begin{equation}}
\def\ee{\end{equation}}
\def\ba{\begin{array}}
\def\ea{\end{array}}
\def\bea{\begin{eqnarray}}
\def\eea{\end{eqnarray}}
\def\beaa{\begin{eqnarray*}}
\def\eeaa{\end{eqnarray*}}
\def\bsp{\begin{split}}
\def\esp{\end{split}}
\def\hh{\!\!\!\!}
\def\EM{\hh &   &\hh}
\def\EQ{\hh & = & \hh}

\def\GE{\hh & \ge & \hh}

\def\GT{\hh & > & \hh}
\def\AND#1{\hh & #1 & \hh}

\def\al{\alpha}

\def\la{\lambda}

\def\ga{\gamma}

\def\da{\delta}

\def\e{\varepsilon}
\def\th{\theta}

\def\ka{\kappa}

\def\oo{\infty}

\def\bfm{{\bf M}}

\def\mcc{{\mathcal C}}
\def\mcl{{\mathcal L}}

\def\mcm{{\mathcal M}}

\def\rd{\,{\rm d}}
\def\dt{\,{\rm d}t}
\def\ds{\,{\rm d}s}

\def\dy{\,{\rm d}y}

\def\dmu{\,{\rm d}\mu}

\def\Lp{{\mathcal L}^p}

\def\Lo{{\mathcal L}^1}
\def\nlo{\|\cdot\|_1}
\def\nlp{\|\cdot\|_p}
\def\nli{\|\cdot\|_\infty}
\def\nmv#1{\|#1\|_{{\bf V}}}

\def\nmo#1{\|#1\|_1}
\def\nmt#1{\|#1\|_2}
\def\nmi#1{\|#1\|_\infty}
\def\nmp#1{\|#1\|_p}
\def\nbv{\|\cdot\|_{\bf V}}
\def\bpr{B_p[r]}
\def\bor{B_1[r]}

\def\sor{S_1[r]}
\def\szr{S_0[r]}
\def\spr{S_p[r]}

\def\ch{\check}
\def\z{\left}
\def\y{\right}
\def\lb{\label}
\def\x#1{(\ref{#1})}

\def\Dx#1{#1^\bullet}
\def\Ddx#1{{\rm d} #1^\bullet}

\def\oo{\infty}

\def\d{\cdot}
\def\dd{\cdots}
\def\pa{\partial}

\def\f{\frac}
\def\q{\quad}
\def\qq{\qquad}
\def\qqf{\quad \forall}
\def\pto{{p\downarrow 1}}

\def\ifl{\iffalse}
\def\nn{\nonumber}
\def\Proof{\noindent{\bf Proof} \quad}
\def\qed{\hfill $\Box$ \smallskip}

\def\ii{\int_I}
\def\ch{\check}

\def\xip{\xi_p}
\def\etap{\eta_p}

\def\ab{\alpha,\beta}

\def\chqr{\check{q}_r}

\def\zzwz#1#2#3#4#5#6#7#8{\bibitem{#1} {#2}, {#3}, {\it #4}, {\bf #5} (#8), {#6}--{#7}.}
\def\shu#1#2#3#4#5#6{\bibitem{#1} {#2}, {\it #3}, {#4}, {#5}, {#6}.}

\def\RS{Riemann-Stieltjes\ }

\def\ysr#1{{\color{red} #1}}

\begin{document}

\title{
Characterization of Maximizers for Sums of the First Two  Eigenvalues of Sturm-Liouville Operators
}

\author{
Gang Meng$^1$, \qq Yuzhou Tian$^{2,}$\footnote{
	Corresponding author.}, \qq Bing Xie$^3$,\qq
Meirong Zhang$^4$
}

\date{}%

\maketitle

\begin{center}
	$^1$ School of Mathematical Sciences, University of Chinese
	Academy of Sciences,\\ Beijing 100049, China\\
	
    $^2$ Department of Mathematics, Jinan University, Guangzhou 510632, China \\

   $^3$ School of Mathematics and Statistics, Shandong University, Weihai, \\ Shandong 264209, China \\
	
	$^4$ Department of Mathematical Sciences, Tsinghua University, Beijing 100084, China\\
	
	E-mail: {\tt
		menggang@ucas.ac.cn (G. Meng), \\ tianyuzhou2016@163.com (Y. Tian),  \\
		xiebing@sdu.edu.cn (B. Xie), \\
	    zhangmr@tsinghua.edu.cn (M. Zhang)
	}
\end{center}

\vspace{-2mm}

\begin{center}
\parbox{14.5cm}{\begin{abstract}
In this paper  we study the maximization of the sum of the first two Dirichlet eigenvalues for Sturm-Liouville operators with potentials in the noncompact space $L^1$. We prove that there exists a unique potential function achieving the maximum, which is non-negative, piecewise smooth, and symmetric.
  Using measure differential equations and weak$^*$ convergence, we show that  the nonzero part of the maximizer can be  determined by  the solution to  the  pendulum equation $\theta'' + \ell \sin\theta = 0 $.

\end{abstract}}

\end{center}

\begin{center}
\parbox{14.5cm}{{\bf Keywords:}
Sturm-Liouville operator, eigenvalue maximization,  measure differential equations, pendulum equation,  nonlinear Schr\"odinger equations}
\end{center}

\begin{center}
\parbox{14.5cm}{{\bf MSC}:
{34L15, 47N10.}}
\end{center}

%

\section{Introduction}
\setcounter{equation}{0} \lb{i}

Given a (real) potential $q$ in the Lebesgue space $\Lp := L^p(I,\R)$, where $I=[0,1]$ and $1\le p < \oo$, we consider the Sturm-Liouville problem
  \be\lb{line}
  y'' + (\la +q(t)) y=0,\qq t\in I,
  \ee
with the Dirichlet boundary condition
  \be \lb{dbc}
  y(0)=y(1)=0.
  \ee
It is well-known that problem \x{line}-\x{dbc} admits a sequence of eigenvalues
  \[
  \la_1(q) < \la_2(q) < \dd < \la_m(q) < \dd, \qq \la_m(q)\to +\oo.
  \]
Since the seminal work by mathematician Krein \cite{Kr55}, Sturm-Liouville problems have remained a central topic in spectral theory and differential equations. Among the numerous works on Sturm-Liouville problems, the main progress consists of  eigenvalue gaps \cite{ MER1991, DM12, GMYZ22, M03, IBSS85,TWZ23,Ze05}, individual eigenvalues \cite{CMZ23,ZW24}, ratios \cite{MR1218744,MR4700368,MR4339006}, nodes \cite{chukm2025,CMWZ24,MR4990941,GZ23} and applications \cite{MR3428972}. 


Besides the eigenvalue problems mentioned above, scholars also focus significantly on eigenvalue sums of differential operators, because they are highly beneficial for understanding the geometry, PDEs and P\'{o}lya conjecture. See \cite{MR4379307,MR2846268,MR4252029,MR3392903} and the references therein.  In a recent paper \cite{TZ25}, they proved that the optimal bounds  for  the sums of the first $m$ eigenvalues of problem \x{line}-\x{dbc} can be attained by some potentials in $\Lp$ with exponent $p>1$. Furthermore, they established a theoretical framework to find these optimal potentials. For exponent $p=1$, the optimal bounds for eigenvalue sums of problem \x{line}-\x{dbc} are essentially different from those for $p>1$, because the $\Lo$ space lacks local compactness. Even the sum of the first two eigenvalues is still unknown. This is the motivation for this paper.

In this paper, we study the maximization problem for the first two eigenvalues
\be \lb{Mr1}
\ch\bfm(r) := \sup\z\{ \la_{1}(q)+\la_2(q): q\in \sor\y\}
=  \sup\z\{ \la_{1}(q)+\la_2(q): q\in \bor\y\},
\ee
where,  for $r\in(0,\oo)$,
\[
\sor :=\z\{ q\in \Lo: \nmo{q} = r\y\},\qq \bor:=\z\{q\in \Lo: \nmo q\le r\y\}
\]
are the sphere and balls in $(\Lo,\nlo)$. Since the $\Lo$ space lacks local compactness, it is a priori unclear whether the maximization problem \x{Mr1} is attainable. If problem \x{Mr1} can be attained by some potentials $\chqr$, it is natural to ask about the uniqueness and symmetry of $\chqr$, as well as its method of construction. Our main result shows that the maximum potential $\chqr$ exists, and it not only has uniqueness and symmetry, but also exhibits a profound connection with the pendulum equation.

\bb{thm}\lb{M21}
For any $r>0$,
there exists a unique potential $\chqr\in \sor$ such that
the maximum in \x{Mr1} is attained by $\chqr$, i.e.
\be\lb{Mr11}
\ch\bfm(r) =\la_1(\chqr)+\la_2(\chqr).
\ee Meanwhile, this potential  $\chqr$ is non-positive, piecewise smooth, and symmetric.
Moreover, there exists a constant $\ch c$ such that
the nonzero part of $\chqr$ is given by
\be \lb{chqr}
\chqr(t) = \ch c + \ell \cos \th(t),
\ee
where the length $\ell$ is
\be \lb{ell}
\ell:=\la_2(\chqr)-\la_1(\chqr)>0,
\ee
and $\th=\th(t)  \in (-\pi, \pi)$ is a segment of a solution to the pendulum
\be \lb{pen}
\th'' + \ell \sin \th=0.
\ee
\end{thm}

The key to proving the above result involves extending the problem to the framework of measure differential equations (MDEs) and leveraging weak$^*$ convergence of measures. By some tricky skills, we analyze the limiting behavior of the critical systems for $\Lp, p>1$ potentials as $\pto$. Furthermore, it is of particular interest that some fundamental equation categories, such as  Ambrosetti-type equations and nonlinear Schr\"odinger equations, play a crucial role in our proof.  We remark that by means of the method developed in this paper, it can be obtained  that the following corresponding minimization problem is unattainable $$\inf\z\{ \la_{1}(q)+\la_2(q): q\in \sor\y\}. $$

This paper is organized as follows. In Section \ref{d}, we present some facts on measures and measure differential equations, including solution regularity and eigenvalue properties. Then we recall some known facts about  critical systems for $\Lp, p>1$ potentials. In Section \ref{cons}, we  analyze the limiting behavior of the critical systems  as $\pto$, and  constructs the maximizing potential $\chqr$ for the case $p = 1$, showing its piecewise smoothness and connection to pendulum equation solutions.

\section{Preliminaries}
\setcounter{equation}{0} \lb{d}

\subsection{Measures and measure differential equations}\lb{mandmde}

In order to describe the oscillation of strings whose distributions may not be absolutely continuous with respect to the Lebesgue measure, one must extend the theory of ODEs to the so-called {\it Generalized Ordinary Differential Equations}\/ \cite{Sc92, S15}, which are involved in the  more general integration theory like the Perron-Stieltjes integral. In his doctoral dissertation, Meng has used only the Riemann-Stieltjes integral \cite{CvB00, DS58} to give a complete extension of ODE to the so-called second order linear {\it measure differential equations}\/ (MDEs). The main results were published in \cite{MZ13}.

Denote by $\mcc:=C(I,\R)$ the Banach space of continuous functions of $I$ with the supremum norm $\nli$. By a measure $\mu$
of $I$, it means that $\mu$ is an element from the dual space $\mcm_0=\mcm_0(I,\R):= (\mcc,\nli)^*$. By the Riesz representation theorem \cite{DS58}, the measure space $\mcm_0$ can be characterized using functions of bounded variations. More precisely,
    \[
    \mcm_0=\bigl\{ \mu:I\to \R: \mu(0+)=0, \ \mu(t+) =\mu(t) \ \forall
t\in (0,1), \ \nmv \mu <\oo\bigr\},
    \]
where $\mu(t+)=\lim_{s\downarrow t}\mu(s)$ is the right-limit and $\nmv \mu$ is the total variation defined by
    \[
    \nmv\mu :=\sup\z\{\sum_{i=0}^{n-1} |\mu(t_{i+1})-\mu(t_i)|:
0=t_0 < t_1< \cdots <t_n=1, \ n\in \N \y\}.
    \]
For example, the unit Dirac measures $\da_a$, $a\in[0,1]$ are as follows. For $a=0$,
    \be \lb{de0}
    \da_0(t):=\z\{\ba{ll} -1 & \mbox{ for } t=0 \\
    0 & \mbox{ for } t\in (0,1] \ea \y. \in \mcm_0.
    \ee
For $a\in (0,1]$,
    \be \lb{dea}
    \da_a(t):=\z\{\ba{ll} 0 & \mbox{ for } t\in [0,a) \\
    1 & \mbox{ for } t\in [a,1] \ea \y. \in \mcm_0.
    \ee
Moreover, $\nmv{\da_a} =1$ and $r \da_a\in \szr$ for all $a\in [0,1]$ and all $r\in [0,\oo)$.

 The second-order scalar linear MDE with a measure $\mu\in \mcm_0$ is written as
    \be \lb{mde}
    \Ddx y + y\dmu(t)=0,\qq t\in I.
    \ee
According to \cite{MZ13}, with an initial value $(y(0),\Dx y(0))=(y_0,z_0)$, the solution $y(t)$ and its generalized right-derivative $\Dx y(t)$ of MDE \x{mde} are determined by the following system of integral equations
    \bea \lb{sd1}
    y(t) \EQ y_0+\int_{[0,t]} \Dx y(s) \ds,\qq t\in [0,1],\\
    \lb{sd2}
    \Dx y(t) \EQ \z\{ \ba{ll} z_0, & t=0,\\
    z_0 - \int_{[0,t]} y(s) \dmu(s), & t\in (0,1],\ea \y.
    \eea
where the integrals are the Lebesgue and the \RS ones respectively \cite{CvB00,DS58}. Both $y$ and $\Dx y{}$ are uniquely defined on $I$.
 Some regularity results on solutions $y(t)$ of MDEs are as follows.

    \bb{rem} \lb{w100}
Any solution $y(t)$ is itself absolutely continuous on $I$ and
    \[
    y'(t) =\Dx y(t) \qq \mbox{a.e. $t\in I$}.
    \]
The generalized velocity $\Dx y(t)$ is the classical right-derivative at all $t\in (0,1)$, i.e.
    \be\lb{dd}
    \Dx y(t) = y'_+(t), \qq t\in(0,1).
    \ee
Moreover, $\Dx y{} \in \mcm:=\mcm(I,\R)$ is a (non-normalized) measure on $I$. Therefore $\Dx y(t)$ is bounded on $I$ and consequently, $y$ belongs to the Sobolev space $W^{1,\oo}(I)$.
    \end{rem}

From the explanation \x{sd1}-\x{sd2} to solutions $y(t)$ of MDE \x{mde}, we have for all $0<t_0 < t\le 1,$
    \bea\lb{sol1}
    y(t) \EQ y(t_0)+\int_{(t_0,t]} \Dx y(s) \ds, \\
    \lb{sol2}
    \Dx y(t) \EQ \Dx y(t_0) - \int_{(t_0,t]} y(s) \dmu(s).
    \eea
Here \x{sol1} is involved in the Lebesgue integrals so that the integrals over $(t_0,t]$ and over $[t_0,t]$ are coincident. However, \x{sol2} is involved in the Lebesgue-Stieltjes integrals so that the integrals over $(t_0,t]$ and over $[t_0,t]$ may be different.

Associated with MDE \x{mde} is the eigenvalue problem
    \be \lb{mula}
    \Ddx y + \la y \dt + y\dmu(t)=0,
\qq t\in I.
    \ee
In \cite{MZ13}, some basic theory for eigenvalues of \x{mula} has been established by using several topological ideas. For example, with the Dirichlet boundary condition \x{dbc}, eigenvalues of problem \x{mula} constitute a (real) sequence  $\{ \la_m(\mu) \}_{m\in \N}$ such that $\la_m(\mu) \nearrow +\oo$ as $m \to +\oo$. In case $\mu\in \mcm_0$ has density $ \mu'=q\in \mcl^1$, solutions of \x{mde} and eigenvalues of \x{mula} return to the corresponding objects of \x{line}.

In the space $\mcm_0$ of measures, there are two topologies. One is induced by the norm $\nbv$ of total variation so that $(\mcm_0,\nbv)$ is a Banach space. Since $\mcm_0=(\mcc,\nli)^*$, one has in $\mcm_0$ the weak$^*$ topology indicated by $w^*$, i.e., $\mu_n\to \mu$ in $(\mcm_0, w^*)$ if and only if
    \[
    \int_I u(t) \dmu_n(t) \to \int_I u(t) \dmu(t)\qq\forall u\in
\mcc.
    \]
For example, for the Dirac measures in \x{de0} and \x{dea}, one has $\da_a\to \da_{a_0}$ in $(\mcm_0, w^*)$ as $a\to a_0$. Some continuity results in \cite{MZ13} on eigenvalues $\la_m(\mu)$ of MDE \x{mula} are as follows.

\bb{lem} \lb{dc} {\rm (\cite{MZ13})} Let\/ $m\in \N$. In the weak$^*$ topology $w^*$ of measures, $\la_m(\mu)$ is continuous in\/ $\mu\in(\mcm_0,w^*)$.
\end{lem}


\subsection{The critical systems of the optimization problem}\lb{cri-eq}

For exponent $1< p < \oo$, we use $\Lp := L^p([0,1],\R)$ to denote the Lebesgue space with the norm $\nlp$. For radius $r>0$, the sphere and ball are respectively
  \[
  \spr :=\z\{ q\in \Lp: \nmp{q} = r\y\}, \qq \bpr:=\z\{q\in \Lp: \nmp q\le r\y\}.
  \]
We consider the following optimization problem in $\Lp$ with a constraint
    \be
    	\lb{Mrp}
    	\ch \bfm_p(r) =  \max (\la_1(q)+\la_2(q))\q \mbox{subject to } \nmp{q}=r.
    	\ee
Based on the following limiting result, the answer to problem \x{Mr1} can be obtained.  For proof, one can refer to \cite{WMZ09, Zh09}.
\bb{lem} \lb{inner}
Given any $r\in(0,\oo)$, one has
\be \lb{limit1}
\lim_{\pto} \ch \bfm_p(r) =\ch \bfm(r).
\ee
\end{lem}    	
    	
In a recent paper \cite{TZ25}, some facts about the  problem \x{Mrp} have been obtained as follows. The problem  \x{Mrp} can be attained by some potentials $q_p=q_{p,r}\in \spr$. By denoting
\[
\xi_p=\la_1(q_p),\qq \eta_p=\la_2(q_p),
\]
there are eigenfunctions $y_p(t)$ and $z_p(t)$ associated with $\xi_p$ and $\eta_p$ respectively such that they satisfy  the following nonlinear ordinary differential systems
$$
\z\{ \ba{l}
y''_p  + \xip y_p  -  \z(y_p^2+z_p^2\y)^{p^*-1} y_p= 0,\\
z''_p + \etap z_p  -  \z(y_p^2+z_p^2\y)^{p^*-1} z_p = 0.
\ea\y.
\eqno(E)
$$
Here $p^*:=p/(p-1)\in(1,\oo)$ is the conjugate exponent of $p$.  These eigenfunctions satisfy
\bea \lb{se3p}
\EM \nmt{y_p} = \nmt{ z_p}, \\
\lb{se4p}
\EM \z\|y_p^2+z_p^2\y\|_{p^*} = r^{p-1}.
\eea
Moreover, the  maximizing potential is given by
\be \lb{qy}
q_p(t) = - \z(y_p^2(t) + z_p^2(t)\y)^{p^*-1}, \qq t\in I.
\ee
 With \x{qy} for $q_p(t)$, system $(E)$ can also be written into linear equations of eigenvalue problems
\bea \lb{ypx}
y''_p +\xip y_p + q_p(t) y_p \EQ 0, \qq t\in I,\\
\lb{zpx}
z''_p + \etap z_p + q_p(t) z_p \EQ 0, \qq t\in I.
\eea
Besides the Dirichlet boundary condition \x{dbc} for $y_p$ and $z_p$, i.e.
    \[
    y_p(0)=y_p(1)=0,\qq z_p(0)=z_p(1)=0,
    \]
$y_p$ and $z_p$ satisfy they have the following nodal properties
    \be \lb{node-p}
    y_p(t)>0\ \forall t\in (0,1),\q z_p(t)>0\ \forall t\in (0,\tau_p),\q z_p(t)<0\ \forall t\in (\tau_p,1).
    \ee
Here $\tau_p\in(0,1)$ is the unique zero of $z_p(t)$ which is called the node of $z_p(t)$ \cite{GZ23}.
Consequently, $y_p(t)$ and $z_p(t)$ are determined by finding positive initial velocities
    \be \lb{abp}
    a_p:= y'_p(0)>0, \q b_p:= z'_p(0)>0.
    \ee
Moreover,  $\ch \bfm_p(r)$ are given by
	\be \lb{bfl2} \ch \bfm_p(r) =\xip+\etap.
	\ee
    Finally, let us notice from \cite[Section 3]{TZ25} that $(E)$ is a Hamiltonian system so that $y_p$ and $z_p$ satisfy
    \be \lb{H}
    (y'_p)^2+ (z'_p)^2+\xi_p y_p^2 + \etap z_p^2 -\f{1}{p^*} \z(y_p^2+ z_p^2\y)^{p^*} \equiv a_p^2+b_p^2.
    \ee

\begin{figure}[ht]
\centering
\includegraphics[width=10cm, height=8cm]{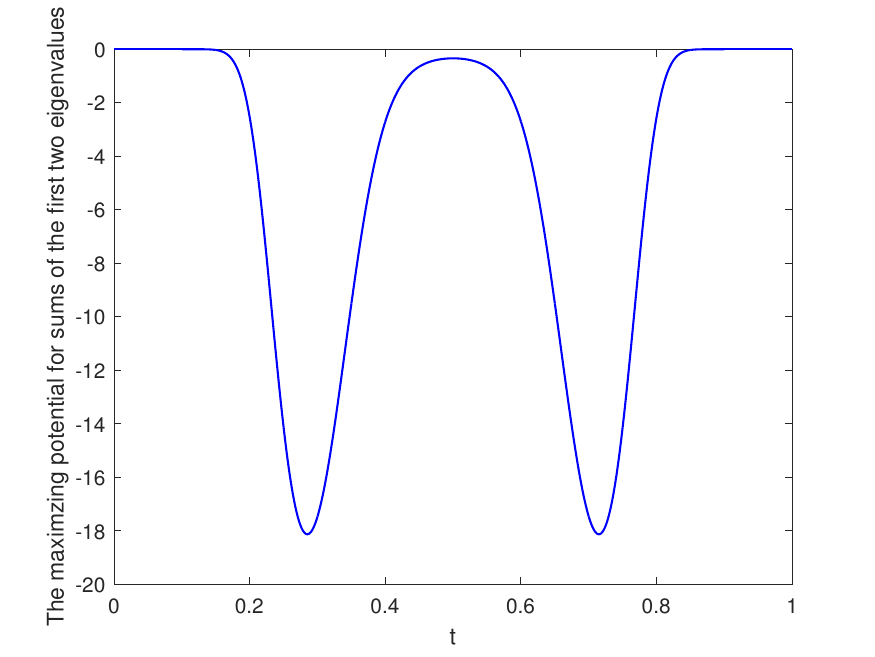}
\caption{The maximizing potential $q_p(t)=q_{p,r}(t)$ with $r=5$ and $p=15/14$.}
\label{maximizerp}
\end{figure}

 We can obtain the maximizing potential $q_p(t)=q_{p,r}(t)$ with fixed $r>0$ and $p>1$  in numerical ways. See Figure \ref{maximizerp} for the case of  $r=5$ and $p=15/14.$

Consider functions
    \bea \lb{ux}
    u_p(t)\AND{:=} y^2_p(t)+z^2_p(t), \qq t\in I,\\
    \lb{kh}
    h_p(t)\AND{:=} 2 a^2_p+ 2 b^2_p-4\xip y^2_p(t) -4\etap z^2_p(t),\qq t\in I.
    \eea
 Trivially, one has
    \be\lb{uh01}
     (u_p(0),u'_p(0))=(u_p(1),u'_p(1)=(0,0), \qq h_p(0)= h_p(1).
    \ee
Moreover,  $u_p(t)>0$ on $(0,1)$. We have the following conclusions.

\bb{lem} \lb{inh1}
Function $u_p(t)$ of \x{ux} satisfies an inhomogeneous linear ODE
    \be \lb{ueq1}
    u''_p + 2(1+1/p^*) q_p(t) u_p = h_p(t),\qq t\in I.
    \ee

\end{lem}

\Proof By using system $(E)$, we have from \x{ux} that
    \beaa
    u'_p \EQ 2 y_p y'_p +2 z_p z'_p, \\
    u''_p \EQ 2 \z( (y'_p)^2 + (z'_p)^2\y) + 2 y_p y''_p+2 z_p z''_p \\
    \EQ 2 \z( a^2_p + b^2_p - \xip y_p^2- \etap z_p^2 + \f{1}{p^*} \z(y_p^2+z_p^2\y)^{p^*}\y)+2 y_p\z( \z(y_p^2+z_p^2\y)^{p^*-1} y_p -\xip y_p\y) \\
    \EM +2 z_p\z( \z(y_p^2+z_p^2\y)^{p^*-1} z_p -\etap z_p\y).
    \eeaa
Here \x{H} and $(E)$ are used. By using \x{ux} and \x{kh}, this is simplified into
    \be \lb{up1}
    u''_p= 2  (1+1/p^*) u_p^{p^*}+ h_p(t).
    \ee
By using the formula \x{qy} for $q_p(t)$, it arrives at Eq. \x{ueq1}. \qed

  \bb{rem}\lb{a-p}
Eq. \x{up1} can also be written as a nonlinear Ambrosetti-type equation
    \be \lb{ueq}
    u''_p - 2 (1+1/p^*) |u_p|^{p^*}= h_p(t),\qq t\in I.
    \ee
 Due to \x{uh01}, both $h_p(t)$ and $u_p(t)$ can be extended to $1$-periodic functions. Hence $u_p(t)$ can be understood as a positive $1$-periodic solution to ODE \x{ueq}.
\end{rem}

\section{Main results}
\setcounter{equation}{0} \lb{cons}

The most important steps of the limiting analysis to system $(E)$ can follow the preceding works by our group. See \cite{WMZ09, Zh09, Zh12}. Let us begin with the following fact.

    \bb{lem}\lb{l-ii}
Let $u_\ga\in \mcc$, $\ga\in(1,\oo)$, be a family of continuous functions on $I$ such that, as $\ga\to+\oo$,
    \be \lb{uc6}
    u_\ga \to u \q\mbox{in } (\mcc,\nli).
    \ee
Then
    \be \lb{uc1}
    \lim_{\ga\to+\oo} \|u_\ga\|_\ga =\|u\|_\oo. 
    \ee
    \end{lem}

\Proof When $u_\ga$ is independent of $\ga$, result \x{uc1} is a well-known fact. For general case as in \x{uc6}, a minor change of the original proof can yield \x{uc1} as well.
\qed

Recall that potentials $q_p(t)$ induce absolutely continuous measures
    \[
    \mu_p(t):= \int_{[0,t]} q_p(s)\ds = -\int_{[0,t]} \z(y^2_p(s)+z^2_p(s)\y)^{p^*-1}\ds, \qq t\in I.
    \]
Since $\nmv{\mu_p}=\nmo{q_p}\le \nmp{q_p} =r$ for all $p\in(1,\oo)$, w.l.o.g., we can assume that, as $\pto$,
     \be \lb{qpw}
    \mu_p \to \mu \qq \mbox{in } (\mcm_0(I),w^*).
    \ee
We will show that $\mu\in \szr$ in Lemma \ref{nm1}  below. By the complete continuity of eigenvalues of MDEs, one has from \x{qpw} that
    \be \lb{llp}
    \xip=\la_1(q_p) \to \la_1(\mu)=: \xi, \qq \etap=\la_2(q_p) \to \la_2(\mu)=: \eta.
    \ee
Here $\xi$ and $\eta$ are the first and the second Dirichlet eigenvalues of the MDE with the limiting measure $\mu$ \cite{MZ13}.

We have the following important observation on $y_p(t)$ and $z_p(t)$ by using the theory for  MDEs \cite{MZ13}. The proof is adopted from \cite{Zh12}.

    \bb{lem} \lb{ucs}
As $\pto$, one has
    \be \lb{ypzp}
    y_p \to y \ne 0,\qq z_p \to z\ne 0 \q \mbox{in } (\mcc,\nli).
    \ee
    \end{lem}

\Proof
Let $\hat y_p(t)$ and $\hat z_p(t)$ be respectively the solutions of Eqs. \x{ypx} and \x{zpx} satisfying the initial values
    \(
    (y(0),y'(0))=(z(0),z'(0))=(0,1).
    \)
Due to \x{qpw}, we have known from the complete continuity result of solutions of MDEs \cite{MZ13} that
    \be \lb{hypzp}
    {\hat y}_p \to {\hat y},\qq {\hat z}_p \to {\hat z} \q \mbox{in } (\mcc,\nli).
    \ee
Here $\hat y(t)$ and $\hat z(t)$ are respectively the solutions of the limiting MDEs
    \bea \lb{y-mu}
    \rd \Dx y + \xi y \dt + y \dmu(t)\EQ 0,\qq t\in I,\\
    \lb{z-mu}
    \rd \Dx z + \eta z \dt + z \dmu(t)\EQ 0, \qq t\in I,
    \eea
with the initial values $(y(0),\Dx y(0))=(z(0),\Dx z(0))=(0,1)$.  In particular, both ${\hat y}(t)$ and ${\hat z}(t)$ are nonzero.

For solutions $y_p(t)$ and $z_p(t)$, by using the initial velocities \x{abp}, we have then from linear ODEs \x{ypx} and \x{zpx} that
    \be\lb{aa}
    y_p(t) \equiv a_p {\hat y}_p(t), \qq z_p(t) \equiv b_p {\hat z}_p(t).
    \ee
By the requirement \x{se3p}, we have from \x{aa} that
    \[
    a_p\nmt{\hat y_p} = \nmt{y_p}=\nmt{z_p} = b_p \nmt{\hat z_p}.
    \]
By \x{hypzp}, we obtain
     \be \lb{abk}
    \f{b_p}{a_p}= \f{\nmt{{\hat y}_p}}{\nmt{{\hat z}_p}}\to \f{\nmt{{\hat y}}}{\nmt{{\hat z}}}=:\ka \in(0,\oo) \qq \mbox{as }\pto.
    \ee
This shows that $b_p$ has the same order as $a_p$ when $\pto$.

We assert that
    \be \lb{abk1}
    \liminf_\pto a_p> 0.
    \ee
Otherwise, if \x{abk1} fails, we will have $a_p\to 0$ and also from \x{abk} that $b_p\to 0$. By \x{hypzp} and \x{aa}, we obtain
    \[
    y_p \to 0,\qq z_p \to 0 \q \mbox{in } (\mcc,\nli).
    \]
To apply Lemma \ref{l-ii} to \x{se4p}, we can choose $\ga=p^*\to+\oo$ and $u_\ga:=y_p^2+ z_p^2 \to 0$ in $(\mcc,\nli)$. The resulted limit will then be
    \[
    \z\|y_p^2+ z_p^2\y\|_{p^*} \to \|0\|_\oo= 0.
    \]
On the other hand, we have from \x{se4p} that $\z\|y_p^2+ z_p^2\y\|_{p^*} \to 1$, a contradiction.

We also assert that
    \be \lb{abk2}
    \limsup_\pto a_p<+\oo.
    \ee
Otherwise, if \x{abk2} fails, we will have $a_p\to +\oo$. Note that
    \[
    \z\|y_p^2+ z_p^2\y\|_{p^*} = a_p^2 \z\|{\hat y}_p^2+ (b_p/a_p)^2 {\hat z}_p^2\y\|_{p^*}.
    \]
By using \x{aa}---\x{abk}, we have
    \[
    \z\|{\hat y}_p^2+ (b_p/a_p)^2 {\hat z}_p^2\y\|_{p^*} \equiv \f{r^{p-1}}{a_p^2}.
    \]
As $\pto$, by applying Lemma \ref{l-ii} again, we obtain
    $$
    \|{\hat y}^2+\ka^2 {\hat z}^2\|_\oo=\lim_\pto \z\|{\hat y}_p^2+ (b_p/a_p)^2 {\hat z}_p^2\y\|_{p^*}= \lim_\pto\f{r^{p-1}}{a_p^2} =0.
    $$
This is impossible, because $\ka {\hat z}(t)$ and ${\hat y}(t)$ are different eigenfunctions for MDE and $\|{\hat y}^2+\ka^2 {\hat z}^2\|_\oo\ne 0$. Hence we have proved \x{abk2}.

Combining \x{abk1} and \x{abk2}, we conclude the following convergence results
    \be \lb{abc9}
    \lim_\pto a_p=a \in(0,+\oo), \qq \lim_\pto b_p=b \in(0,+\oo).
    \ee
Finally, we can obtain from \x{hypzp}, \x{aa} and \x{abc9} that
    \be \lb{yzp}
    y_p = a_p {\hat y}_p\to a {\hat y}=:y, \qq z_p = b_p {\hat z}_p\to b {\hat z}=:z
    \qq \mbox{in }(\mcc,\nli).
    \ee
These are the uniform convergence results stated in \x{ypzp}.
\qed

    \bb{rem} \lb{la9}
In the proof of Lemma \ref{ucs}, $y_p$ and $z_p$ are considered as solutions of initial value problems of ODEs \x{ypx} and \x{zpx}. Due to \x{yzp}, $y$ and $z$ are solutions of MDEs \x{y-mu} and \x{z-mu} with the initial values $(y(0),\Dx y (0))=(0,a)$ and $(z(0), \Dx z(0))= (0,b)$ respectively. On the other hand, $y_p$ and $z_p$ are the Dirichlet eigenfunctions of ODEs, the limiting solutions $y$ and $z$ are respectively eigenfunctions associated with eigenvalues $\la_1(\mu)=\xi$ and $\la_2(\mu)=\eta$. Therefore  $y(t)$ and $z(t)$ have similar nodal properties as in \x{node-p}.
    \end{rem}

    \bb{lem}\lb{nm1}
The limiting eigenfunctions $y(t)$ and $z(t)$ in Lemma \ref{ucs} and the limiting measure $\mu$ in \x{qpw} satisfy
    \bea \lb{le1}
    \nmt{y}\EQ \nmt{z},\\
    \lb{le2}
    \nmi{y^2+z^2} \EQ 1,\\
    \lb{le3}
    \ii (y^2+z^2) \dmu \EQ -r.
    \eea
In particular, one has $\mu\in \szr$.
    \end{lem}

\Proof Due to \x{ypzp}, equality \x{le1} comes from \x{se3p} immediately.

Applying Lemma \ref{l-ii} to \x{se4p} with the choice $\ga=p^*\to \oo$, we obtain \x{le2}.

By using \x{qy}, equality \x{se4p} is
    \[
    - r^p = \ii -\z(y_p^2(t)+z_p^2(t)\y)^{p^*} \dt= \ii q_p(t)\d (y_p^2(t)+z_p^2(t)) \dt.
    \]
By letting $\pto$, the left-hand side has the limit $ -r$. Due to \x{ypzp}, one has $y_p^2+z_p^2 \to y^2+z^2$ in $(\mcc,\nli)$. Combining with the weak$^*$ convergence \x{qpw}, we know that the right-hand side has the limit $\ii (y^2+z^2) \dmu$. Hence one has result \x{le3}.

It is well-known from the weak$^*$ convergence \x{qpw} that the variation norm of $\mu$ satisfies
    \[
    \nmv{\mu} \le \liminf_{\pto} \nmv{\mu_p}\le r.
    \]
Moreover, we have from \x{le2} and \x{le3} that
    \[
    r=-  \ii (y^2+z^2) \dmu \le \nmi{y^2+z^2} \nmv{\mu} = \nmv{\mu}.
    \]
These have shown that $\nmv{\mu}=r$, i.e. $\mu\in \szr$. \qed

The following limiting inhomogeneous linear MDE \x{mde3} is one of the important tools of this paper.

\bb{lem}\lb{umu} As $\pto$,  one has
    \bea \lb{ut}
    u_p \AND{\to} u:=y^2+z^2,\\
    \lb{thz}
     h_p\AND{\to} h:=  2a^2+2b^2-4\xi y^2 -4\eta z^2
    \eea
in $(\mcc,\nli)$. Moreover, the limiting function $u$ satisfies the following inhomogeneous linear MDE
    \be \lb{mde3}
    \rd \Dx u +  2 u\dmu(t) = h(t) \dt,\qq t\in I,
    \ee
and the initial values $(u(0),\Dx u(0))=(0,0)$.
\end{lem}

\Proof As $\pto$, we have from \x{ypzp} the result \x{ut}. Moreover, it follows from \x{abc9} that
    \[
    a^2_p + b^2_p\to a^2+b^2.
    \]
Then we obtain from \x{kh} that
    \(
    h_p\to 2a^2+2b^2-4\xi y^2 -4\eta z^2
    \)
in $(\mcc,\nli)$. This proves \x{thz}.

Let us consider linear ODE \x{ueq1} as an inhomogeneous linear MDE.  By Lemma 3.9 in \cite{MZ13},  $u_p$ is the solution to problem  \x{ueq1} with $u_p(0) =0, u_p'(0) =0$ iff it
	satisfies, for each $ t\in [0,1],$
	\bea \lb{fix90}
	\EM  u_p(t) = \int_0^1 G(t,s) ( -2(1+1/p^*) q_p(s) u_p(s) +  h_p(s)) \ds \nn
	\\ \EQ   -\int_0^1 G(t,s)  2(1+1/p^*) u_p(s) \dmu_p(s) + \int_0^1 G(t,s)  h_p(s) \ds,
	\eea
	where   $G: [0,1] \times [0,1] \to \R$ is defined as
	\[
	G(t,s) := \z\{ \ba{ll}t-s & \mbox{ for } 0 \le s \le t \le 1, \\
	0 & \mbox{ for } 0 \le t < s \le 1. \ea \y.
	\]
	Because of the convergence results \x{qpw}, \x{ut} and \x{thz}, we have from \x{fix90} that, as $p \to 1,$ \beaa
	u(t) =  -\int_0^1 G(t,s)  2 u(s) \dmu(s) + \int_0^1 G(t,s)  h(s) \ds \q \forall t\in [0,1].
	\eeaa
	Thus, the limiting MDE \x{mde3} follows. \qed

We notice from \x{le2} and \x{ut} that $\nmi{u}=\nmi{y^2+z^2}=1$. Let us introduce the following subsets of $I$
    \[
    I_0: = \z\{ t\in I: y^2(t)+z^2(t)< 1\y\}, \qq
    I_1: = \z\{ t\in I: y^2(t)+z^2(t)= 1\y\}.
    \]
Then both $I_0$ and $I_1$ are non-empty, and $I$ has the disjoint decomposition
    \(
    I= I_0 \cup I_1 .
    \)
The set $I_0$ is a (relatively) open subset of $I$, $\{0,\, 1\}\subset I_0$, and $I_0$ is composed of at most countably many open subintervals.
Meanwhile, the set $I_1\subset (0,1)$ is closed, compact, and $I_1$ is composed of at most countably many closed subintervals which may shrink into single points. In particular, $\pa I_1\ne \emptyset$ is a countable subset.
Moreover, one has the following properties.

\bb{lem}\lb{cst}
The limiting measure $\mu$ satisfies ${\rm supp}(\mu) \subset I_1$. Thus $\mu$ is massless on $I_0$ and  $\mu(t)$ is constant on each open subinterval $I_0$. Moreover, the solutions $y(t)$ and $z(t)$ and the measure $\mu(t)$ are $C^\oo$ on $I_0$.
\end{lem}

\Proof Consider any maximal subinterval $\ch I$ of $I_0$. For $t\in \ch I$, one has $y^2(t)+z^2(t)<1$ and $q_p(t) = -\z(y_p^2(t) +z_p^2(t)\y)^{p^*-1} \to 0$. Because the limiting measure $\mu$ is non-increasing, the set of its continuous points is dense in $\ch I$. For any continuous points $t_1,t_2\in \ch I$, we have by the Alexandroff Theorem  \cite[p. 316]{DS58} that \beaa  \mu(t_2) - \mu(t_1) = \int_{[t_1,t_2]} \dmu  = \lim_\pto   \int_{[t_1,t_2]} \dmu_p = \lim_\pto  \int_{t_1}^{t_2} q_p(s)\ds = 0,\eeaa
which implies that $\mu$ is constant on $\ch I$  .   Hence the measure $\mu$ has the density
    \be \lb{q-0}
    q(t):=\mu'(t)=0,\qq t\in I_0.
    \ee
Because of \x{q-0}, MDEs \x{y-mu} and \x{z-mu} are then reduced to ODEs
    \bea \lb{y-0}
    y'' +\xi y\EQ 0,\qq t\in I_0,\\
    \lb{z-0}
    \qq z''+\eta z\EQ 0,\qq t\in I_0.
    \eea
From \x{q-0}, \x{y-0} and \x{z-0}, we know that $y(t)$, $z(t)$ and $\mu(t)$ are $C^\oo$ on $I_0$. \qed

\bb{lem}\lb{smty}
Suppose that $I_1$ has a nontrivial maximal subinterval $[\ab]$. Then the solutions $y(t)$ and $z(t)$ and the measure $\mu(t)$ are $C^\oo$ on $(\ab)$. Moreover, on $(\ab)$, $y(t)$ satisfies ODE
 \be \lb{E1}
y'' + \z(\xi+a^2+b^2 -2\eta \y )  y +2 (\eta-\xi) y^3=0,
\ee
and the measure $\mu(t)$ has the density given by
    \be \lb{q-y}
  q(t):= \mu'(t) = a^2+b^2 -2\eta +2 (\eta - \xi) y^2.
  \ee
\end{lem}

\Proof In the sequel, we consider $t\in(\ab)$. By \x{dd}, one has
    \[
    \Dx u(t) = u'_+(t) =0.
    \]
Now we use the explanations to solutions of MDE \x{mde3} like in \x{sd1} and \x{sd2}. See also \x{sol1} and \x{sol2}. Hence  \beaa
0\EQ \Dx u(t) = \Dx u(\al) - 2 \int_{(\al,t]} u(s)\dmu(s) +\int_{(\al,t]} h(s) \ds\\
\EQ - 2 \int_{(\al,t]}\dmu(s) +\int_{[\al,t]} h(s) \ds\\
\EQ - 2 \mu(t) + 2 \mu(\al) +\int_{[\al,t]} h(s) \ds.
\eeaa
Since $h(t)$ is $C^1$ on $I$, we know that $\mu(t)$ is $C^2$ on $(\ab)$. Consequently, we have from \x{thz} that
        \[\bb{split}
    q(t)&=\mu'(t) = \f{1}{2} h(t)=   a^2+b^2-2\xi y^2 -2\eta z^2   \\
    & = a^2+b^2-2\xi y^2 -2\eta (1-y^2) = a^2+b^2 -2\eta +2 (\eta - \xi) y^2,
    \end{split}
    \]
because $y^2(t)+z^2(t)=1$ on $[\ab]$. This gives relation \x{q-y}.
  Moreover, with the help of \x{q-y}, MDE \x{y-mu} becomes the ODE
        \[
       y'' + \xi y = - y q(t) = -y\z( a^2+b^2 -2\eta +2 (\eta - \xi) y^2 \y).
       \]
This gives ODE \x{E1}. From the theory of ODEs, we know that, as a solution to ODE \x{E1}, $y(t)$ is $C^\oo$ on $(\ab)$.  The $C^\oo$ smoothness of $\mu(t)$  follows from \x{q-y}. Finally,  $z(t)$ is also $C^\oo$, because $z(t)$ satisfies the eigenvalue equation
    \[
     z'' +\eta z +q(t) z = 0  , \qq t\in(\ab).
    \] \qed

The equation  \x{E1} is a one-dimensional cubic nonlinear Schr\"odinger equation. This type of equation is also used in \cite{WMZ09, Zh09} in determining the optimal bounds for eigenvalues of Sturm-Liouville operators.

Since $\pa I_1$ is at most countable, we can use \x{q-0} and \x{q-y} to define an integrable potential $\ch q$ on $I$. Such a potential $\ch q$ is piecewise $0$ or $C^\oo$ on subintervals of $I$ and  defines an absolutely continuous measure
    \be \lb{muac}
    \mu_{\rm ac}(t)=\int_{[0,t]} \ch q(s)\ds\in \mcm_0.
    \ee
At boundary points from $\pa I_1$, the measure $\mu(t)$ may admit nonzero Dirac measures. Hence we have the following decomposition for the measure $\mu$.

\bb{lem} \lb{mu} There exists a possible sequence $\{\tau_i\}\subset \pa I_1$ of different positions and a possible sequence $\{r_i\}$ of nonzero masses which define a completely singular measure
    \be \lb{mucs}
    \mu_{\rm cs}=\sum_i - r_i \da_{\tau_i}.
    \ee
Moreover, combining with the absolutely continuous measure in \x{muac}, we know that $\mu$ is decomposed into
    \be \lb{mus}
    \mu=\mu_{\rm ac} + \mu_{\rm cs}.
    \ee
Furthermore, one has
    \bea  \lb{taui}
    y^2(\tau_i)+z^2(\tau_i) \EQ 1\qqf i,\\
    \lb{ris}
     r_i \GT 0\qqf i,\\
    \lb{in90}
    - \ch q(t) = |\ch q(t)|\GE 0\qqf t\in I,\\
    \lb{mubvr}
    \nmv{\mu} = \nmo{\ch q} +\sum_i |r_i| \EQ r.
    \eea
\end{lem}

\Proof We have already explained the decomposition results \x{muac}---\x{mus}. Fact \x{taui} comes from $\{\tau_i\}\subset \pa I_1$.

For facts \x{ris} and \x{in90}, let us notice that $\la_1(\nu)+\la_2(\nu)$ is strictly decreasing in measures $\nu \in \mcm_0$. Since $\mu$ has no Dirac measures at endpoints $0$ and $1$, we know that, as a function, for the maximal problem, $\mu(t)$ must be non-increasing in $t\in I$. This implies \x{ris} and \x{in90}.

We have from \x{mus} that
    \[
    \nmv{\mu} =\nmo{\ch q} +\sum_i |-r_i| =\nmo{\ch q} +\sum_i r_i.
    \]
Hence \x{mubvr} can be obtained from the fact $\nmv{\mu}=r$ in Lemma \ref{nm1}. \qed

Now we can obtain an important Lemma about the  maximizer $\ch q$.

\bb{lem} \lb{max1}
For the maximization problem \x{Mr1}, the completely singular part $\mu_{\rm cs}=0$ in the decomposition \x{mus}. That is, problem \x{Mr1} is attained by a potential $\ch q$ defined before. Moreover, $\ch q(t) \le 0$ for all $t\in I$ and $\nmo{\ch q}=r$, i.e. $\ch q\in S_1[r]$.
\end{lem}

\Proof For the maximization problem \x{Mr1}, suppose that one has some negative  mass $ -r_i<0$ in \x{mucs}. We know from the explanations to solutions of MDEs \x{y-mu} and \x{z-mu} that
    \[
    y'(\tau_i+) - y'(\tau_i-) = r_i y(\tau_i), \qq z'(\tau_i+) - z'(\tau_i-) = r_i z(\tau_i).
    \]
Thus
    \beaa
    u'(\tau_i+) -u'(\tau_i-) \EQ \z(2 y(\tau_i) y'(\tau_i+)+ 2 z(\tau_i) z'(\tau_i+) \y) - \z(2 y(\tau_i) y'(\tau_i-) + 2 z(\tau_i) z'(\tau_i-)\y) \\
    \EQ 2 y(\tau_i) \z(y'(\tau_i+)-y'(\tau_i-)\y)+2 z(\tau_i) \z(z'(\tau_i+)-z'(\tau_i-)\y)  \\
    \EQ 2 r_i \z(y^2(\tau_i)+z^2(\tau_i)\y) \\
    \EQ 2 r_i >0.
    \eeaa
This is impossible because $u(t)$ attains its maximum $1$ at $t=\tau_i$.  \qed

Due to Lemma \ref{max1}, $I_1$ must admit some nontrivial subinterval $[\ab]$. Note that the nonlinear Schr\"odinger equation \x{E1} is integrable. Hence the potential $\ch q$ can be determined by \x{q-y} using the solution $y$ of \x{E1}.

Now we are in the position to prove the main results about maximizers for sums of the first two Dirichlet eigenvalues.


\noindent {\bf Proof of Theorem \ref{M21}.} For the uniqueness, let $q_1$ and $q_2$ be two maximizing potentials.  Then, their average $q_3 = \f{1}{2}q_1+ \f{1}{2}q_2 \in \bor.$   Let us introduce the first normalized  eigenfunctions $u_1,u_2, u_3$ associated with the first eigenvalues $\la_1(q_1),\la_1(q_2),\la_1(q_3)$  and the second  normalized eigenfunctions $v_1,v_2,  v_3$ associated with the second eigenvalues $ \la_2(q_1), \la_2(q_2), \la_2(q_3).$ As shown in \cite{Fan1949},
 \beaa\lb{lam} \lambda_1(q)+ \lambda_2(q) \EQ  \min_{u,v \in H^1_0, \langle u,v \rangle =0 } R(u,v,q)
  \\ \EQ    \min_{u,v \in H^1_0, \langle u,v \rangle =0 }  \z( \frac{ \int_0^1 u'^2 \dt - \int_0^1 q u^2 \dt }{ \int_0^1 u^2\dt} + \frac{ \int_0^1 v'^2 \dt - \int_0^1 q v^2 \dt }{ \int_0^1 v^2\dt} \y) ,
\eeaa
where $\langle u,v \rangle = \int_0^1 u(t)v(t)\dt.$
Unless $span\{u_1, v_1\} = span\{u_3, v_3\} $  and $span\{u_2, v_2\} = span\{u_3, v_3\},$  we have
\beaa  \EM \bfm(r) \ge  \la_1(q_3) + \la_2(q_3) =  R(u_3,v_3, q_3)
\\ \EQ  \f{1}{2} R(u_3,v_3,q_1)  + \f{1}{2}     R(u_3,v_3,q_2)
\\ \GT	\f{1}{2} \z(\la_1(q_1) + \ \la_2(q_1) \y)
+ \f{1}{2}  \z( \la_1(q_2) +  \la_2(q_2) \y)
\\ \EQ  \bfm(r).\eeaa
Therefore, $span\{u_1, v_1\} = span\{u_2, v_2\} $  and then $$ q_1(t)  + \la_1(q_1) +  q_1(t)  + \la_2(q_1)= q_2(t)   + \la_1(q_2) + q_2(t)   + \la_2(q_2) \mbox{ a.e. $t\in [0,1]$},$$
which implies that  $q_1(t) = q_2(t)  \mbox{ a.e. $t\in [0,1]$}.$

We claim that the maximizer $\ch q_r(t)$ is symmetric with respect to $1/2$, i.e., $\ch q_r(t) = \ch q_r(1 -t)$ for all $t \in [0,1].$ In fact,  it can be checked that  $ q^*(t) := \ch q_r(1-t)$  also serves as a maximizing potential for the optimization problem \x{Mr1}. Consequently,  the symmetry of the maximizer can be deduced directly from the uniqueness of the maximizer.

Let  $[\ab]$  be any nontrivial maximal subinterval of $I_1$.
 For $t\in(\ab)$, we have from \x{y-mu} and \x{z-mu} the following ODEs
\be \lb{yzeq}
y''+ \xi y + \ch q(t) y=0, \qq z''+\eta z +\ch q(t) z=0.
\ee
By eliminating $\ch q(t)$ from \x{yzeq}, we obtain
\be \lb{yz22}
\f{z''}{z} -\f{y''}{y}+\eta-\xi=0, \qq t\in (\ab)\setminus\{\tau\}.
\ee
Note from $y^2+z^2=1$ that
\[
z = \pm \sqrt{1-y^2}, \qq
z'' = \pm \f{(y^3-y)y''-y'^2}{\z(\sqrt{1-y^2}\y)^3}.
\]
Then Eq. \x{yz22} simplifies to the following second-order nonlinear ODE
\be \lb{E2}
(1-y^2) y''+y y'^2 - (\eta-\xi) y (1-y^2)^2=0
\ee
for $t\in (\ab)\setminus\{\tau\}$. In case $\tau \in(\ab)$, one has $y(\tau)=1$, $y'(\tau)=0$ and therefore \x{E2} still holds at $t=\tau$. Thus Eq. \x{E2} is well-defined for all $t\in(\ab)$.

Note that $y^2(t)+z^2(t)=1$ on $[\ab]$. By Lemma \ref{smty}, $y(t)$ and $z(t)$ are $C^\oo$ on $(\ab)$. Hence there exists a $C^\oo$ function $\th=\th(t)$ such that
\be \lb{po}
(y(t),z(t)) = (\cos(\th(t)/2), \sin(\th(t)/2)), \qq t\in(\ab).
\ee
As $y(t)$ is positive, one can choose   $\th(t) \in (-\pi, \pi)$. Moreover, for $t\in (\ab)\setminus \{\tau\}$, one has from \x{po}
\[
y'' = -\f{1}{2} \th'' \sin(\th/2) - \f{1}{4}\th'^2 \cos(\th/2).
\]
Thus \x{E2} is
\[\bb{split}
0 & = (1-y^2) y''+y y'^2 - (\eta-\xi) y (1-y^2)^2\\
& = -\f{1}{2} (\th''+(\eta-\xi)\sin \th)\sin^3(\th/2).
\end{split}
\]
By discarding the nonzero factor $-\f{1}{2}\sin^3(\th/2)$, we know that $\th(t) \in (-\pi, \pi)$ satisfies the pendulum equation  \be \lb{P}
\th''+(\eta-\xi)\sin \th=0,
\ee for $t\in (\ab)\setminus \{\tau\}$. Since $\th(t)$ is $C^\oo$, Eq. \x{P} is also verified at $t=\tau$.

As for $\ch q(t)$, we obtain from  \x{q-y}
\[\bb{split}
\ch q(t) & =  a^2+b^2 -2\eta +2 (\eta - \xi)  \cos^2(\th(t)/2)\\
& =  a^2+b^2 -\eta -\xi  +(\eta-\xi) \cos\th(t)
\\ &  = \ch c + \ell \cos \th(t),
\end{split}
\] with $\ch c =  a^2+b^2 -\eta -\xi$ and $\ell = \eta-\xi.$
This gives \x{chqr}.\qed

\begin{rem} The maximizing potential and the corresponding eigenfunctions for the case of $r = 5$ can be numerically simulated, as illustrated in Figure \ref{r5}. Similarly, those for $r = 20$ are presented in Figure \ref{r20}. Based on numerical simulations, we conjecture that there exists a constant $r^* (\approx 15.5292)$ such that  the maximizer $\check{q}_{r}$ contains  two nonzero parts when $r<r^*$ and one nonzero part when $r>r^*.$ See Figure \ref{rstar} for the case of $r = r^*.$

\begin{figure}[H]
	\centering
	\includegraphics[width=6.5cm, height=6cm]{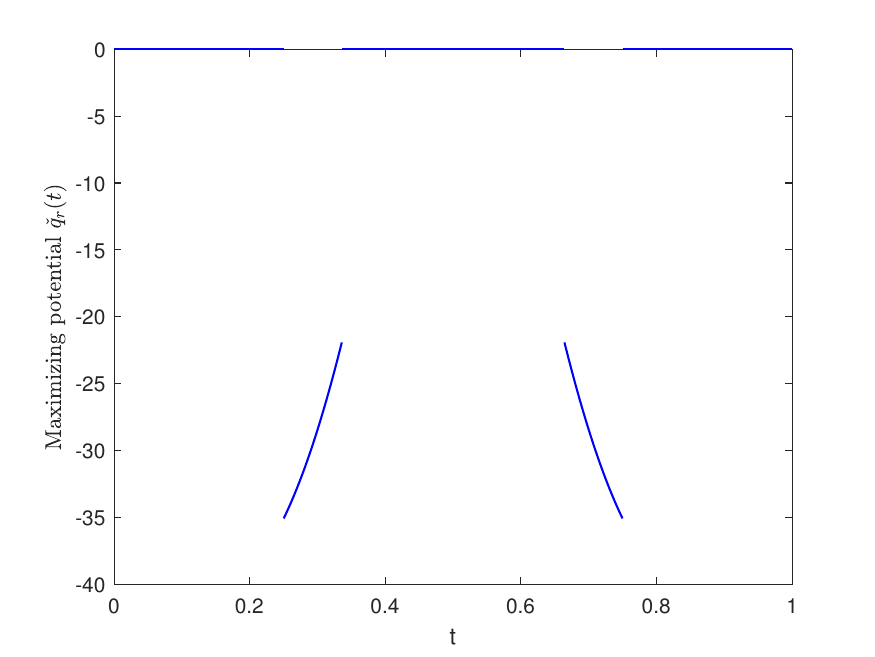}\hspace{1.5cm}
	\includegraphics[width=6.5cm, height=6cm]{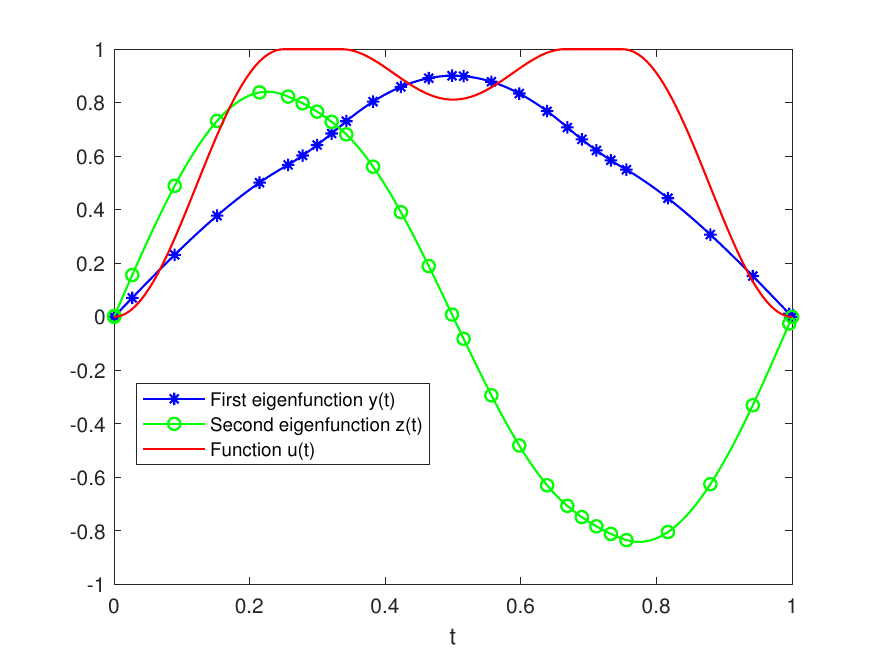}
	\caption{The maximizing potential $\check{q}_{r}(t)$  (left) and  eigenfunctions (right) with $r=5.$
	}
	\label{r5}
\end{figure}

\begin{figure}[H]
	\centering
	\includegraphics[width=6.5cm, height=6cm]{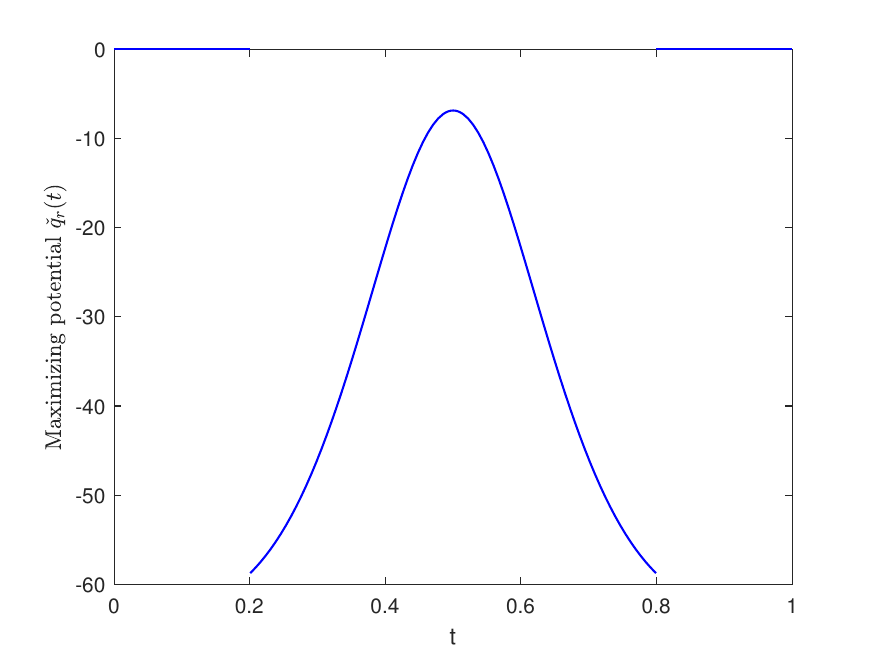}\hspace{1.5cm}
	\includegraphics[width=6.5cm, height=6cm]{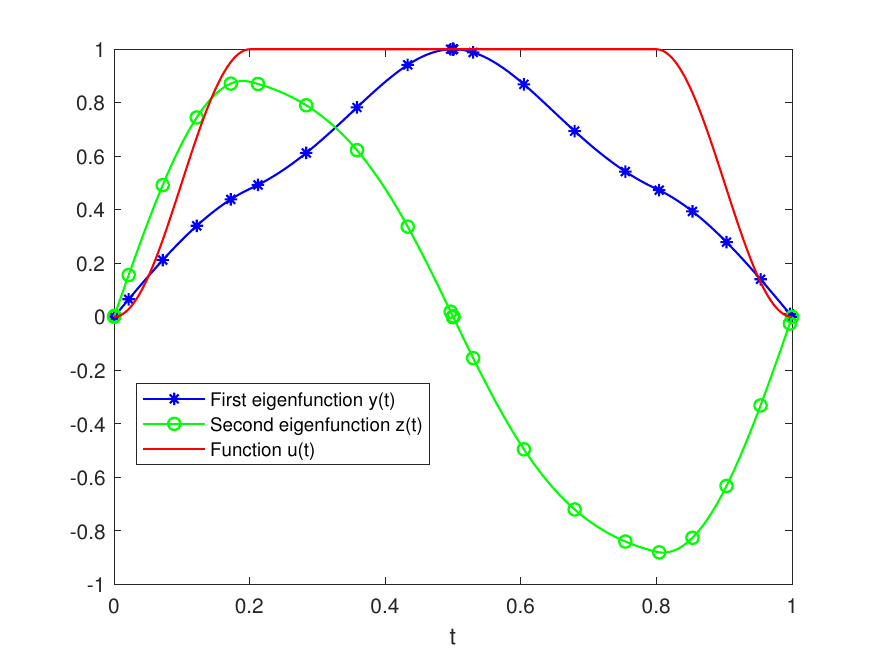}
	\caption{The maximizing potential $\check{q}_{r}(t)$  (left) and  eigenfunctions (right) with $r=20.$
	}
	\label{r20}
\end{figure}

	\begin{figure}[H]
	\centering
	\includegraphics[width=6.5cm, height=6cm]{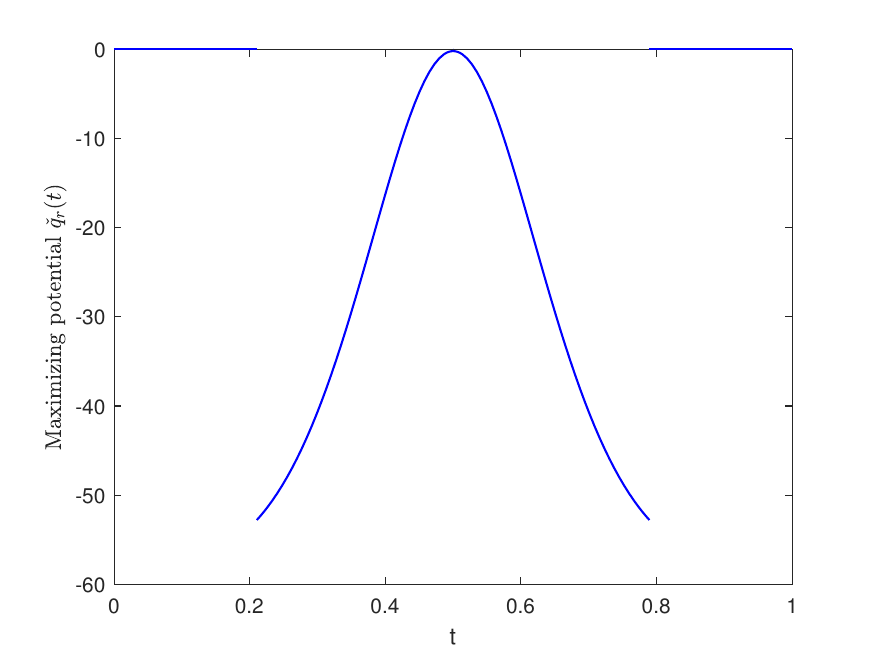}\hspace{1.5cm}
	\includegraphics[width=6.5cm, height=6cm]{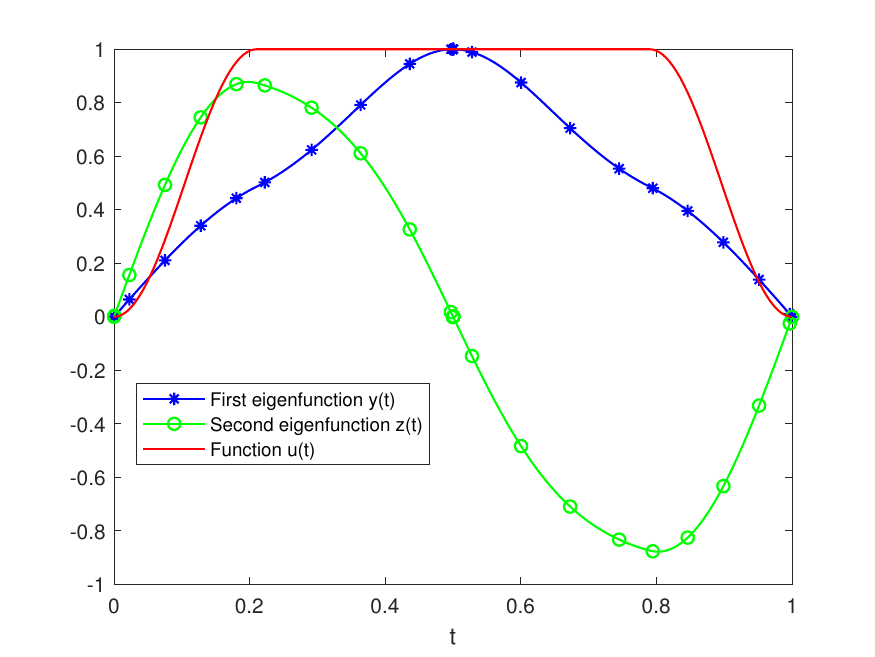}
	\caption{The maximizing potential $\check{q}_{r}(t)$  (left) and  eigenfunctions (right) with $r=r^*.$
	}
	\label{rstar}
\end{figure}

\end{rem}

\section*{Acknowledgments}

This research is supported by  the National Natural
Science Foundation of China (Grant No. 12471174, 12401205 and 12371182), the Guangdong Basic and Applied Basic Research Foundation (Grant no. 2023A1515110430), the Science and Technology Planning Project of Guangzhou (Grant no. 2025A04J3404), the China Postdoctoral Science Foundation (Grant no. 2025M773097), the Postdoctoral Fellowship Program (Grade C) of China Postdoctoral Science Foundation (Grant no. GZC20230970) and the Fundamental Research Funds for the Central Universities (No. 21624347).

\end{document}